\documentclass[a4paper,12pt]{amsart}
\usepackage{amssymb}
\usepackage{amsmath}
\usepackage{ifthen}
\usepackage{graphicx}
\usepackage{amsfonts}
\usepackage{amscd}
\usepackage{amsxtra}
\usepackage{color}

\newtheorem{thm}{Theorem}
\newtheorem{cor}{Corollary}
\newtheorem{lem}{Lemma}

\newtheorem{conj}{Conjecture}
\theoremstyle{definition}
\newtheorem{example}[equation]{Example}
\newtheorem{prob}[equation]{Problem}

\newcounter {own}
\def\theown {\thesection       .\arabic{own}}

\newenvironment{rem}{%
\bigskip
\noindent \textsl{{\sl Remark. }}}{\bigskip}
\newenvironment{rems}{%
\bigskip
\noindent \textsl{{\sl Remarks. }}}{\bigskip}

\newenvironment{pf}[1][]{%
 \vskip 3mm
 \noindent
 \ifthenelse{\equal{#1}{}}%
  {{\slshape Proof. }}%
  {{\slshape #1.} }%
 }%
{\qed\bigskip}

\newcounter{alphabet}
\newcounter{tmp}
\newenvironment{Thm}[1][]{\refstepcounter{alphabet}%
\bigskip%
\noindent%
{\bf Theorem \Alph{alphabet}}%
\ifthenelse{\equal{#1}{}}{}{ (#1)}%
{\bf .} \itshape}{\vskip 8pt}

\makeatletter
\newcommand{\Ref}[1]{\@ifundefined{r@#1}{}{\setcounter{tmp}{\ref{#1}}\Alph{tmp}}}
\makeatother

\newenvironment{Lem}[1][]{\refstepcounter{alphabet}%
\bigskip%
\noindent%
{\bf Lemma \Alph{alphabet}}%
{\bf .} \itshape}{\vskip 8pt}

\newcommand{\IN}{{\mathbb N}}
\newcommand{\IC}{{\mathbb C}}
\newcommand{\ID}{{\mathbb D}}

\newcommand{\real}{{\operatorname{Re}\,}}

\def\be{\begin{equation}}
\def\ee{\end{equation}}

\newcommand{\bee}{\begin{enumerate}}
\newcommand{\eee}{\end{enumerate}}

\newcommand{\blem}{\begin{lem}}
\newcommand{\elem}{\end{lem}}
\newcommand{\bthm}{\begin{thm}}
\newcommand{\ethm}{\end{thm}}
\newcommand{\bcor}{\begin{cor}}
\newcommand{\ecor}{\end{cor}}
\newcommand{\beg}{\begin{example}}
\newcommand{\eeg}{\end{example}}
\newcommand{\begs}{\begin{examples}}
\newcommand{\eegs}{\end{examples}}
\newcommand{\bdefe}{\begin{defin}}
\newcommand{\edefe}{\end{defin}}
\newcommand{\bprob}{\begin{prob}}
\newcommand{\eprob}{\end{prob}}
\newcommand{\bei}{\begin{itemize}}
\newcommand{\eei}{\end{itemize}}

\newcommand{\bcon}{\begin{conj}}
\newcommand{\econ}{\end{conj}}
\newcommand{\bcons}{\begin{conjs}}
\newcommand{\econs}{\end{conjs}}
\newcommand{\bprop}{\begin{propo}}
\newcommand{\eprop}{\end{propo}}
\newcommand{\br}{\begin{rem}}
\newcommand{\er}{\end{rem}}
\newcommand{\brs}{\begin{rems}}
\newcommand{\ers}{\end{rems}}
\newcommand{\bo}{\begin{obser}}
\newcommand{\eo}{\end{obser}}
\newcommand{\bos}{\begin{obsers}}
\newcommand{\eos}{\end{obsers}}
\newcommand{\bpf}{\begin{pf}}
\newcommand{\epf}{\end{pf}}
\newcommand{\ba}{\begin{array}}
\newcommand{\ea}{\end{array}}
\newcommand{\beq}{\begin{eqnarray}}
\newcommand{\beqq}{\begin{eqnarray*}}
\newcommand{\eeq}{\end{eqnarray}}
\newcommand{\eeqq}{\end{eqnarray*}}

\newcommand{\ra}{\rightarrow}

\newcounter{minutes}
\divide\time by 60
\newcounter{hours}
\multiply\time by 60 \addtocounter{minutes}{-\time}

\begin{document}
\bibliographystyle{amsplain}
\title[Extreme points method and univalent harmonic mappings]{Extreme points method and univalent harmonic mappings}

\thanks{
File:~\jobname .tex,
          printed: \number\day-\number\month-\number\year,
          \thehours.\ifnum\theminutes<10{0}\fi\theminutes}

\author[Y. Abu Muhanna]{Yusuf Abu Muhanna}
\address{Y. Abu Muhanna, Department of Mathematics,
American University of Sharjah, UAE-26666.}
\email{ymuhanna@aus.edu}

\author[S. Ponnusamy]{Saminathan Ponnusamy $^\dagger $
}
\address{S. Ponnusamy,
Indian Statistical Institute (ISI), Chennai Centre, SETS (Society
for Electronic Transactions and security), MGR Knowledge City, CIT
Campus, Taramani, Chennai 600 113, India. }
\email{samy@isichennai.res.in, samy@iitm.ac.in}

\keywords{Analytic, univalent, convex, close-to-convex, starlike and
harmonic mappings, functions of bounded boundary rotation, integral means, extreme points. \\
$
^\dagger$ {\tt This author is on leave from the Indian Institute of Technology Madras, India}
}

\begin{abstract}
We consider the class of all sense-preserving complex-valued harmonic mappings $f=h+\overline{g}$
defined on the unit disk $\ID$ with the normalization $h(0)=h'(0)-1=0$ and $g(0)=g'(0)=0$ with the
second complex dilatation $\omega:\,\ID\rightarrow \ID$, $g'(z)=\omega (z)h'(z)$. In this paper, the
authors determine sufficient conditions on $h$ and $\omega$ that would imply the univalence of
harmonic mappings $f=h+\overline{g}$ on $\ID$.
\end{abstract}

\maketitle \pagestyle{myheadings}
\markboth{Y. Abu Muhanna, and S. Ponnusamy}{Extreme points method and univalent harmonic mappings}

\section{Preliminaries and Main Results}
Denote by ${\mathcal A}$ the class of all functions $h$ analytic in the unit disk ${\mathbb D}=\{z \in {\mathbb C}:\, |z|<1\}$
with the normalization $h(0)=0=h'(0)-1$. We let $\mathcal S$ denote the subset of
functions from ${\mathcal A}$ that are univalent in $\ID$.
A locally univalent function $h\in {\mathcal A}$ is in ${\mathcal S}^*(\alpha)$ if and only
if $ {\rm Re } \left(zh'(z)/h(z)\right) > \alpha$ for $z\in \ID$ where $\alpha <1$. A locally univalent function
$h\in {\mathcal A}$ is said to belong to ${\mathcal K}(\alpha)$ if and only if $zh' \in {\mathcal S}^*(\alpha)$. Functions in ${\mathcal S}^*(0)$
and ${\mathcal K}(0)$ are referred to as the normalized starlike (with respect to $h(0)=0$) and convex functions in $\ID$, respectively.
Finally, $h\in {\mathcal A}$ is called close-to-convex if there exists a $g\in {\mathcal S}^*(0)$ such that
$ {\rm Re } \left(e^{i\gamma}zh'(z)/g(z)\right) > 0$ for $z\in \ID$ and for some $|\gamma |<\pi/2$.
It is known that every close-to-convex function is univalent (see \cite{Dur,Pomm}).

We remind the reader that a univalent analytic or harmonic function $f$ on $\ID$ is close-to-convex if $f(\ID)$ is close-to-convex, i.e.
its complement in $\IC$ is the union of closed half lines with pairwise disjoint interiors. We direct the reader
to \cite{Clunie-Small-84,Du} and expository notes \cite{SaRa2013} for several basic knowledge on planar univalent harmonic mappings
and methods of constructing them.

Our first problem concerns the class $ \mathcal{K}(\beta)$ of functions $h\in {\mathcal S}$ such that
\be\label{class-G}
{\rm Re } \left( 1+\frac{zh''(z)}{h'(z)}\right) > \beta ~\mbox{ for $z\in \ID$},
\ee
for some $\beta \in [-1/2,1)$. For convenience, we let $ {\mathcal K}(-1/2)=\mathcal{F}$.
In particular, functions in $\mathcal{F}$ are known to be close-to-convex but are not
necessarily starlike in $\ID$. For $\beta \geq 0$, functions in $\mathcal{K}(\beta)$ are known to be convex
in $\ID$.

%


There are two important sufficient conditions for close-to-convexity of harmonic
mappings due to Clunie and Sheil-Small \cite{Clunie-Small-84}. We now recall them here for a ready
reference.

\begin{Lem}\label{LemA}
If a harmonic mapping $f=h+\overline{g}$ satisfies the condition
$|g'(0)|<|h'(0)|$ and that the analytic function $h+\epsilon g$ is close-to-convex
for each $\epsilon$ ($|\epsilon|=1$), then $f$ is also close-to-convex.
\end{Lem}

\begin{Lem}\label{LemB}
Let $h$ be analytic and convex in $\ID$. If $g$ and $\omega$ are analytic in $\ID$ such that
$|\omega (z)|<1$ and $ g'(z)=\omega (z)h'(z)$ for $z\in\ID$, then
every harmonic mapping of the form $f=h+\overline{g}$ is close-to-convex and univalent in $\ID$.
\end{Lem}

In our proof, we observe that Lemma \Ref{LemB} is an immediate consequence of Lemma \Ref{LemA} (see also the proof of the
case $\beta =0$ of Theorem \Ref{BL-Theorem}(2)).

The function $\omega:\,\ID\rightarrow \ID$ satisfying the relation $ g'= \omega h'$ is called the
second complex dilatation of the sense-preserving harmonic mapping $f=h+\overline{g}$ in $\ID$. In our discussion,
it is convenient to consider harmonic functions $f=h+\overline{g}$ in $\ID$ with the standard normalization, namely,
$h(0)=0=h'(0)-1$ and $g(0)=0$, and the family of normalized harmonic convex (resp. close-to-convex and starlike)
mappings, i.e. sense-preserving univalent harmonic functions that have a convex (resp. close-to-convex and starlike) range
(see \cite{Clunie-Small-84,Du,SaRa2013}).

As an application of Lemma \Ref{LemA} and Kaplan's characterization of close-to-convex functions,
the following results were obtained in \cite{Bshouty-Lyzzaik-2010,BshoutyJoshiJoshi-2013} (see also
Bharanedhar and Ponnusamy \cite{Bhara-samy-pre11}).

\begin{Thm}\label{BL-Theorem}
Let $f=h+\overline{g}$ be a harmonic mapping in $\ID$ such that $ g'(z)=\omega (z)h'(z)$ in $\ID$ for some $\omega :\,\ID\rightarrow \ID$.
Then $f$ is close-to-convex in $\ID$ if one of the following conditions is satisfied:
\bee
\item[{\rm (1)}] $h \in {\mathcal K}(-1/2)$ and $\omega (z)=e^{i\theta}z$ in $\ID$
\item[{\rm (2)}] $h \in {\mathcal K}(\beta)$ for some $\beta \in (-1/2, 0]$ and $|\omega (z)|<\cos (\beta \pi)$ for $z\in\ID$.
\eee
\end{Thm}

Originally, Theorem \Ref{BL-Theorem}(1) was a conjecture of Mocanu \cite{Mocanu-2010} and was settled by
Bshouty and Lyzzaik \cite{Bshouty-Lyzzaik-2010} (see also \cite{Bhara-samy-pre11})
whereas Theorem \Ref{BL-Theorem}(2) extends Lemma \Ref{LemB} (see \cite[Theorem 4]{BshoutyJoshiJoshi-2013}
and \cite{PonSai-14a}). We remark that the case $\beta =0$ of Theorem \Ref{BL-Theorem}(2) is equivalent to Lemma \Ref{LemB}.

In Section \ref{sec-2}, using extreme points method, we present an elegant proof of Theorem \Ref{BL-Theorem}
and several other new results. Second consequence of our method gives for example the following.

\bthm \label{AP-Theorem}
Let $h\in \mathcal{F}$. 
Then for $\beta >0$ and $r\in (0,1)$ one has
$$I_\beta (r,f)=\frac{1}{2\pi}\int_{-\pi}^{\pi}\frac{d\theta}{|h'(re^{i\theta})|^{2\beta }}
\leq \frac{2^{6\beta }}{\pi}{\bf B}\Big (\frac{6\beta +1}{2},\frac{1}{2}\Big ),
$$
where ${\bf B}(.,.)$ denotes the Euler-beta function. The inequality is sharp.
\ethm

We present two different proofs of Theorem \ref{AP-Theorem}. One of the proofs relies
on the method of extreme points (see for example \cite{HM}) while the other relies on the subordination relation.
The estimates of $I_1$ has received special attention in the field of planar fluid mechanics,
where these functionals are participating in isoperimetric problems for moving phase domains, eg. \cite{Vasi2002}
and \cite{VasiMark2003}.

In order to present the third consequence of our approach, we recall the following result
which is a partial extension of the classical result of Alexander's theorem from conformal mappings to univalent harmonic mappings.

\begin{Thm}\label{Thm-Alex1}
{\rm (\cite[p.108, Lemma]{Du})}
Let $f=h+\overline{g}$ be a sense-preserving harmonic starlike mapping in $\ID$. If $H$ and $G$ are the analytic functions
defined by the relations
\be\label{ps4eq11}
zH'(z) = h(z),~ ~zG'(z) = - g(z),~ ~H(0)=G(0)=0,
\ee
then $F=H+\overline{G}$ is a convex mapping in $\ID$.
\end{Thm}

A generalization of Theorem \Ref{Thm-Alex1} has been obtained by Ponnusamy and Sairam Kaliraj \cite{PonSai-14a}. However, it is natural
to ask what would be the conclusion if the assumption about $f$ is replaced just by the analytic part $h$ being starlike in $\ID$.
We remark that the harmonic Koebe function (see \cite{Clunie-Small-84,Du,SaRa2013}) $K$ defined by
$$K(z)=\frac{z-\frac{1}{2}z^2+\frac{1}{6}z^3}{(1-z)^3}+
\overline{\left (\frac{\frac{1}{2}z^2+\frac{1}{6}z^3}{(1-z)^3}\right )} \quad \mbox{ for $z\in \ID$},
$$
is starlike in $\ID$ whereas its analytic part is not even univalent in $\ID$. Also, there are harmonic convex function
whose analytic part is not necessarily starlike in $\ID$.

\bthm\label{Thm-Alex2}
Let $f=h+\overline{g}$ be a sense-preserving harmonic mapping in $\ID$, where $h\in {\mathcal S}^*$ and $g(0)=0$. If $H$ and $G$ are
the analytic functions defined by the relations \eqref{ps4eq11}, then for each $|\lambda|\leq 1$,
the harmonic function $F_{\lambda} =H+\lambda\overline{G}$ is sense-preserving and close-to-convex mapping in $\ID$. In particular,
$F =H+\overline{G}$ is a close-to-convex mapping in $\ID$.
\ethm

We now state our next result whose proof follows similarly. So we omit its detail.

\bthm\label{Thm-Alex3}
Let $f=h+\overline{g}$ be a harmonic mapping in $\ID$, where $h\in {\mathcal S}^*(\beta)$ for some
$\beta \in (-1/2,0]$, $g(0)=0$ and $ g'(z)=\omega (z)h'(z)$ in $\ID$ for some $\omega :\,\ID\rightarrow \ID$ satisfying the condition
$|\omega (z)|<\cos (\beta \pi)$ for $z\in\ID$. If $H$ and $G$ are
the analytic functions defined by the relations \eqref{ps4eq11}, then for each $|\lambda|= 1$,
the harmonic function $F_{\lambda} =H+\lambda\overline{G}$ is sense-preserving and close-to-convex mapping in $\ID$.
\ethm

We remark that functions in ${\mathcal S}^*(\beta)$ are not necessarily univalent in $\ID$ if $\beta <0$.
At the end of the article, Bshouty and Lyzzaik \cite{Bshouty-Lyzzaik-2010} expressed their
interest in determining sufficient condition on $h$ so that $ g'(z)=e^{i\theta}zh'(z)$ implies that $f=h+\overline{g}$ is univalent in $\ID$.
Several of the remaining results of this article motivate their desire by choosing $h$ appropriately.
Proof of Theorem \ref{Thm-Alex2} will be given in Section \ref{sec-2}.

%
%
%



For $\alpha \in (1,2)$, let $ CO_{H}(\alpha)$ denote the class of all
harmonic mappings $f=h+\overline{g}$ defined on $\ID,$ where
$ g'(z)= \omega (z)h'(z)$ with $|\omega (z)|<1$ for $z\in \ID$ and
$h\in CO(\alpha )$, the class of all concave univalent functions (see Section \ref{sec-2} for the precise
definition).
The class $ CO(\alpha )$ has been extensively studied in the recent years and for a detailed discussion
about concave functions, we refer to \cite{Avk-Wir-06, Avk-Wir-05, BPW-09, Pom-Cruz} and the references
therein. We now state our next result.

\bthm\label{AP-Theorem3}
For $\alpha \in (1,2)$, let $f=h+\overline{g}\in CO_{H}(\alpha)$. If the dilatation $\omega $ satisfies the conditions
$|\omega (z)|< \sin (\frac{2-\alpha }{2})\pi$ for $z\in\ID$, then $f$ is close-to-convex (univalent) in $\ID$.
\ethm

A simple consequence of Theorem \ref{AP-Theorem3} gives

\bcor\label{AP-Cor3}
For $\alpha \in (1,2)$, each harmonic mapping $f\in CO_{H}(\alpha)$ with the dilatation
$\omega (z)=\left( \sin (\frac{2-\alpha }{2})\pi \right) e^{i\theta}z$ for $z\in\ID$ is close-to-convex (univalent) in $\ID$.
\ecor

We conjecture that Corollary \ref{AP-Cor3} is sharp in the sense that the number $ \sin (\frac{2-\alpha }{2})\pi$
cannot be replaced by a larger one for a given $\alpha \in (1,2)$.

A function $h$ analytic and locally univalent in $\ID$
is said to have {\it boundary rotation} bounded by $K\pi$, $K\ge 2$, if for $0<r<1$
\begin{equation}\label{sec2-eq1}
\int_0^{2\pi}\left|\real \left (1+\frac{re^{i\theta}h''(re^{i\theta})}{h'(re^{i\theta})}
\right )\right|\,d\theta \le K\pi\,.
\end{equation}
Let ${\mathcal V}_K$ be the class of all analytic functions $h$ in $\ID$ (with the normalization
$h(0)=0=h'(0)-1$) having boundary rotation bounded by $K\pi$.
The reader is referred to Paatero \cite{Paa1} (see also \cite{Dur,Koepf-89} and Section \ref{sec-2} for additional information
about the class ${\mathcal V}_K$), where the study of these classes was initiated, for the geometric significance.

\bthm\label{AP-Theorem4}
Let $h\in {\mathcal V}_{K}$ with $2\leq K\leq 4-\delta $ for a fixed $\delta \in [0,2]$, and the dilatation satisfies the
condition $|\omega (z)|<\sin (\frac{\delta \pi }{4})$ for $z\in\ID$. Then harmonic mapping $f=h+\overline{g}$ is
close-to-convex and univalent in $\ID$.
\ethm

As an immediate corollary to this result, we have

\bcor \label{AP-Cor2}
Let $f=h+\overline{g}$ be a harmonic mapping in $\ID$ such that $h\in {\mathcal V}_{K}$ with $2\leq K\leq 4-\delta $ for a fixed $\delta \in [0,2]$,
and that $g'(z)=e^{i\theta }\sin (\frac{\delta \pi }{4})z h'(z)$ for $z\in\ID$. Then $f$ is close-to-convex and univalent in $\ID$.
\ecor

We conjecture that Corollary \ref{AP-Cor2} is sharp in the sense that the number $\sin (\frac{\delta \pi }{4})$ cannot be replaced by
a larger one for a given $\delta <2$.

\begin{figure}
\begin{center}
\includegraphics[height=6.0cm, width=5.5cm, scale=1]{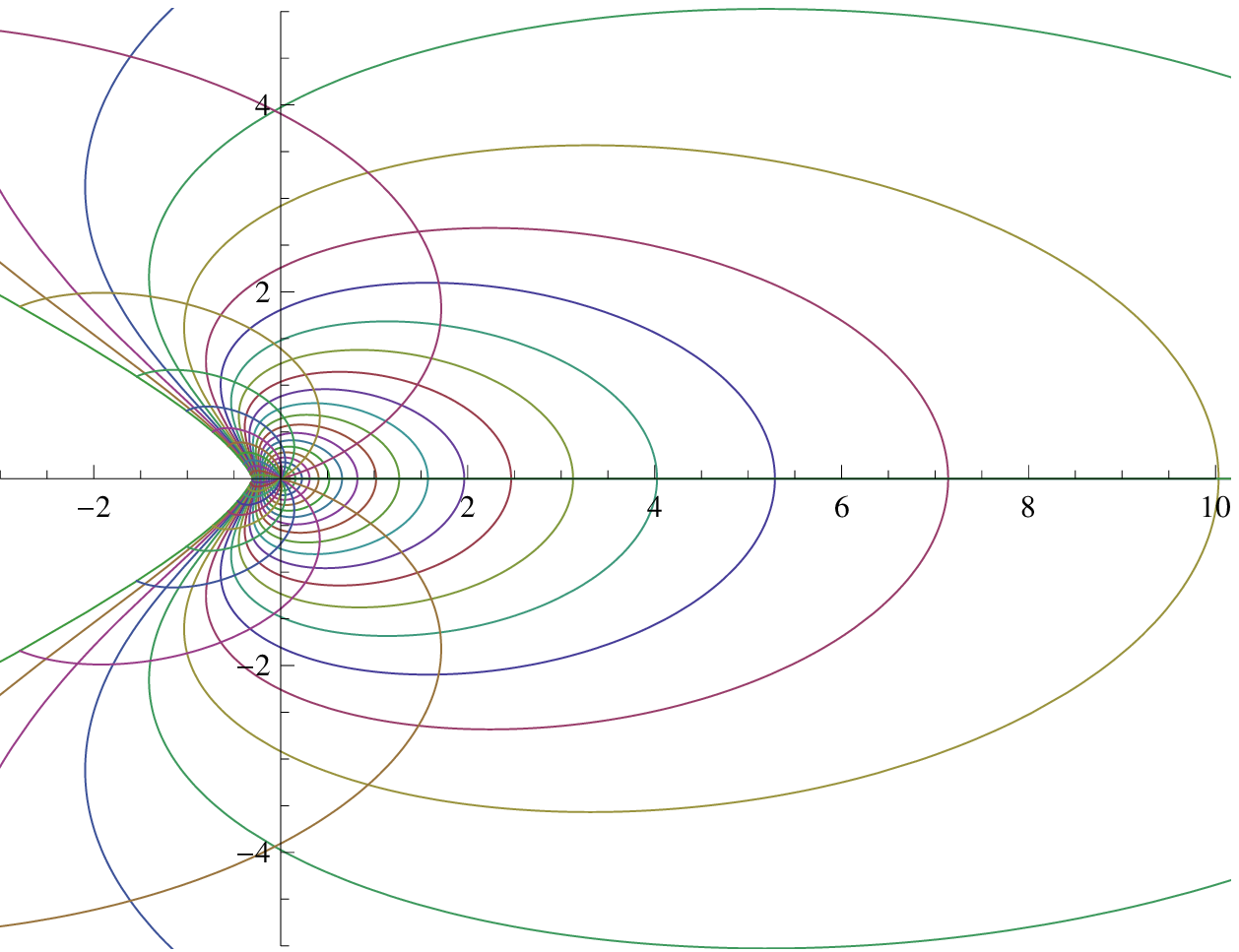}
\hspace{1cm}
\includegraphics[height=6.0cm, width=5.5cm, scale=1]{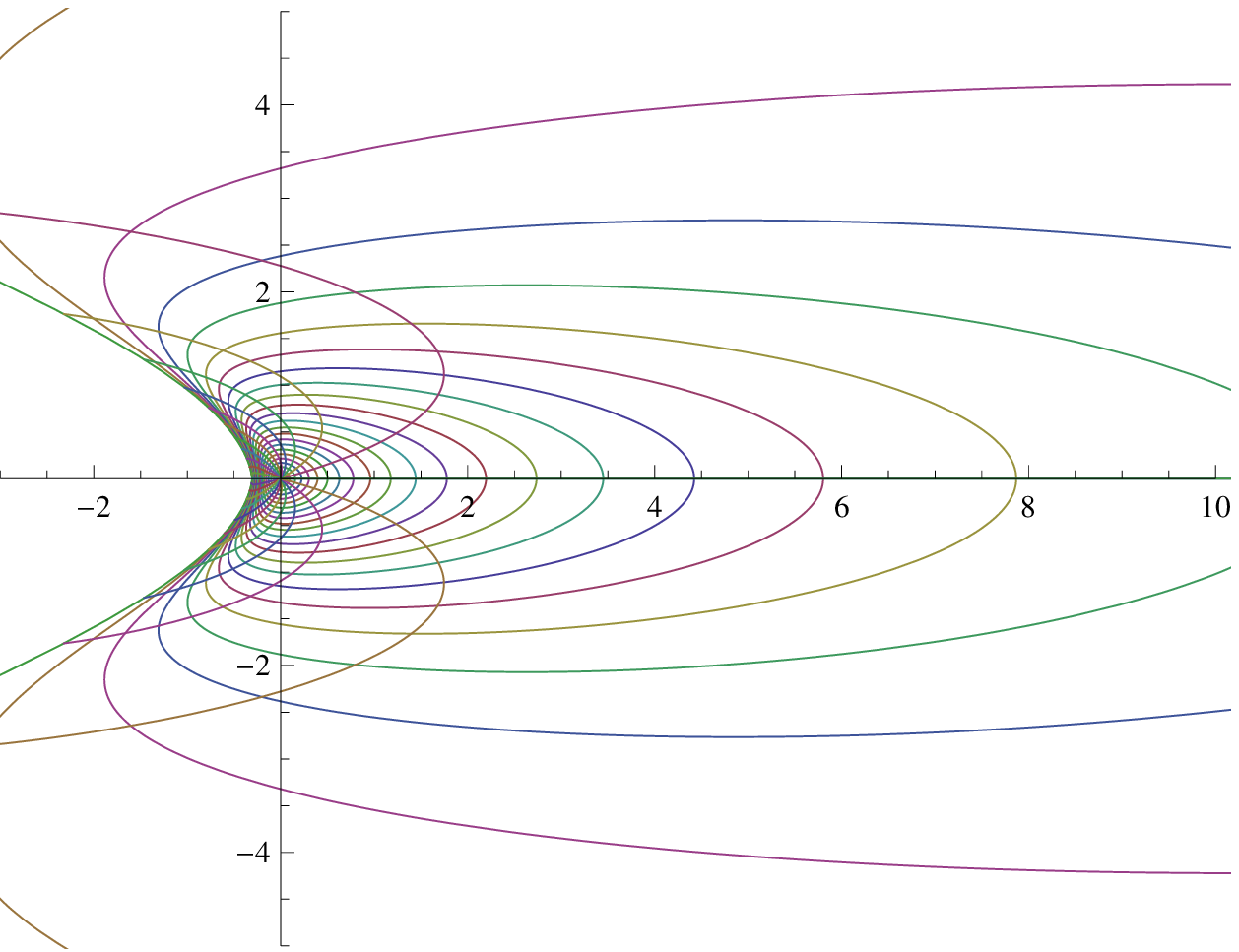}
\end{center}
(a) $\delta = 0.5$, $K = 3$ \hspace{5cm} (b) $\delta = 1 $, $K= 2.5$
\vspace{0.5cm}

\begin{center}
\includegraphics[height=6.0cm, width=5.5cm, scale=1]{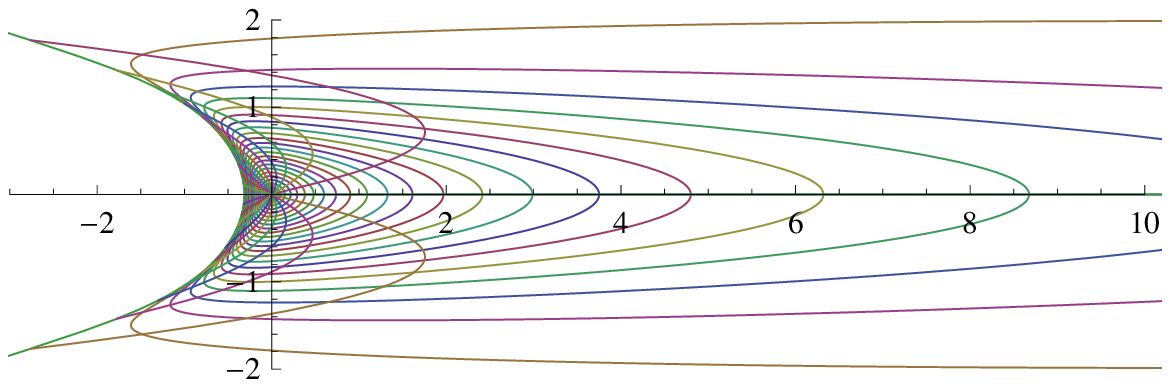}
\hspace{1cm}
\includegraphics[height=6.0cm, width=5.5cm, scale=1]{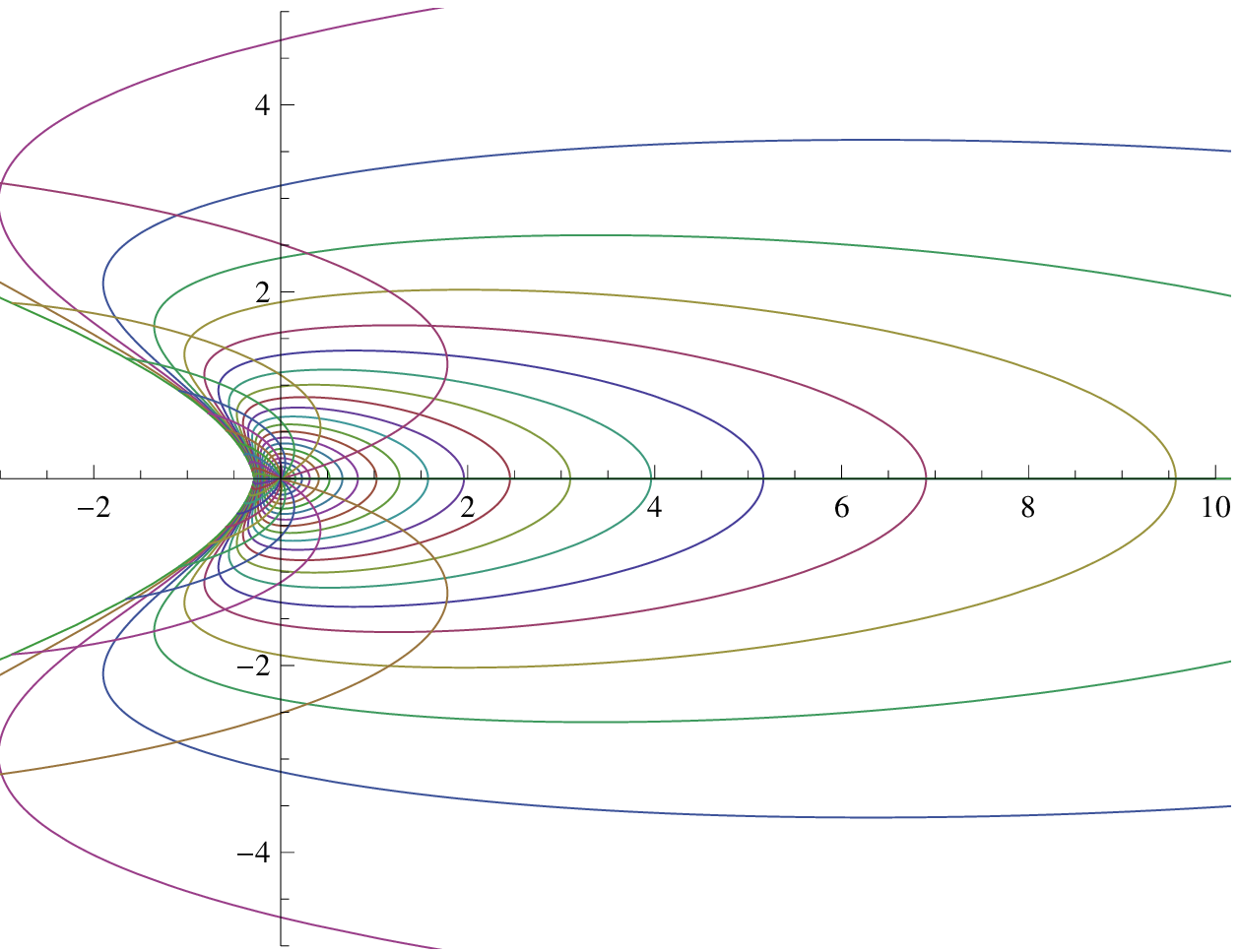}
\end{center}
(c) $\delta = 1.5$, $K =2.1 $ \hspace{5cm} (d) $\delta = 1$, $K =2.75 $

%

\begin{center}
\includegraphics[height=6.0cm, width=5.5cm, scale=1]{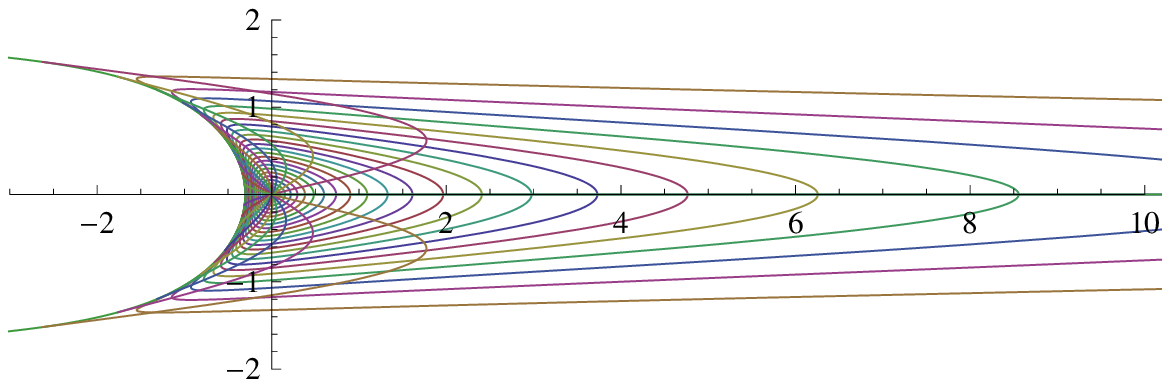}
\hspace{1cm}
\includegraphics[height=6.0cm, width=5.5cm, scale=1]{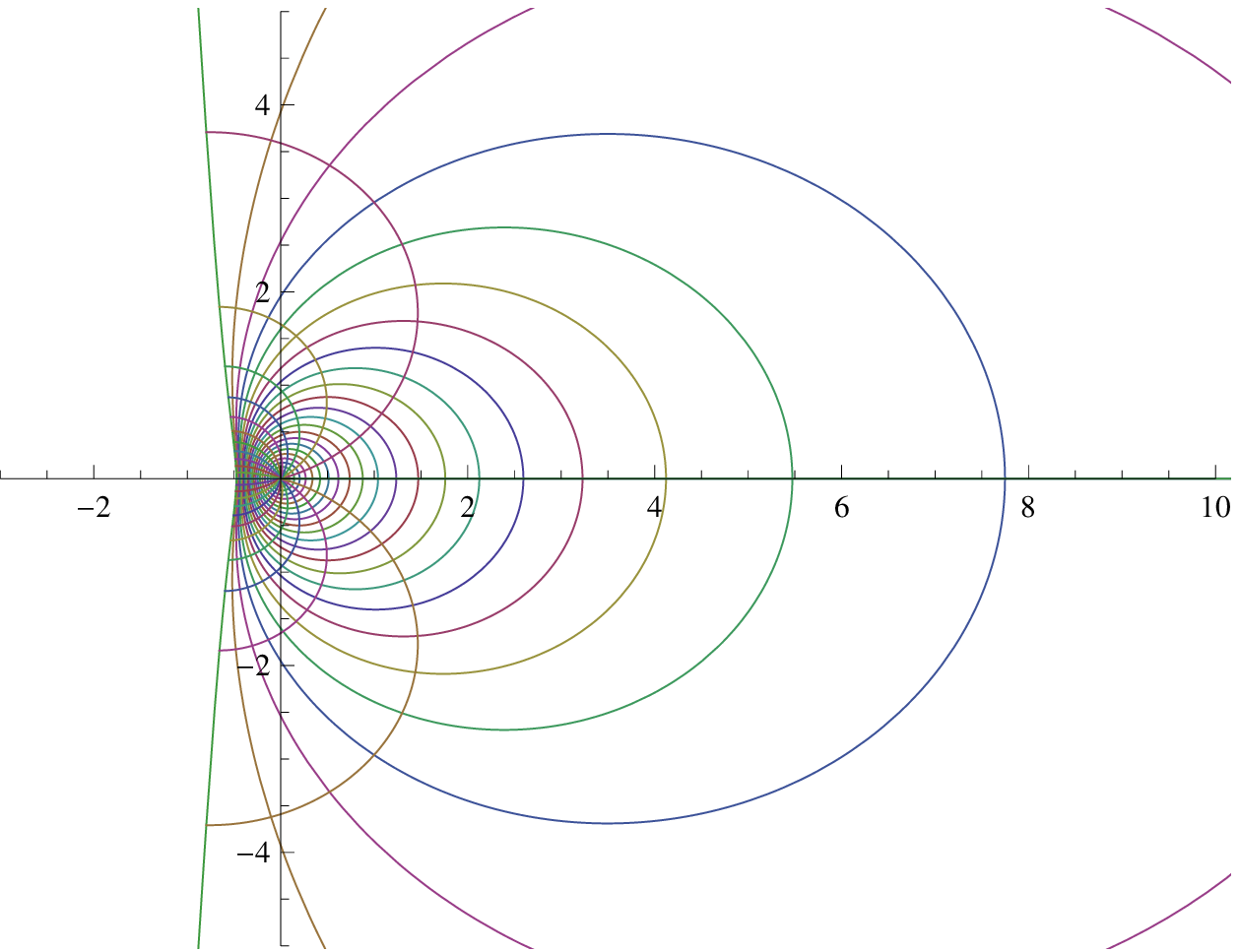}
\end{center}
(e) $\delta = 1.9$, $K = 2.05$ \hspace{5cm} (f) $\delta = 0.1 $, $K = 2.05 $
\vspace{0.1cm}
\end{figure}

\begin{figure}
\begin{center}
\includegraphics[height=6.0cm, width=5.5cm, scale=1]{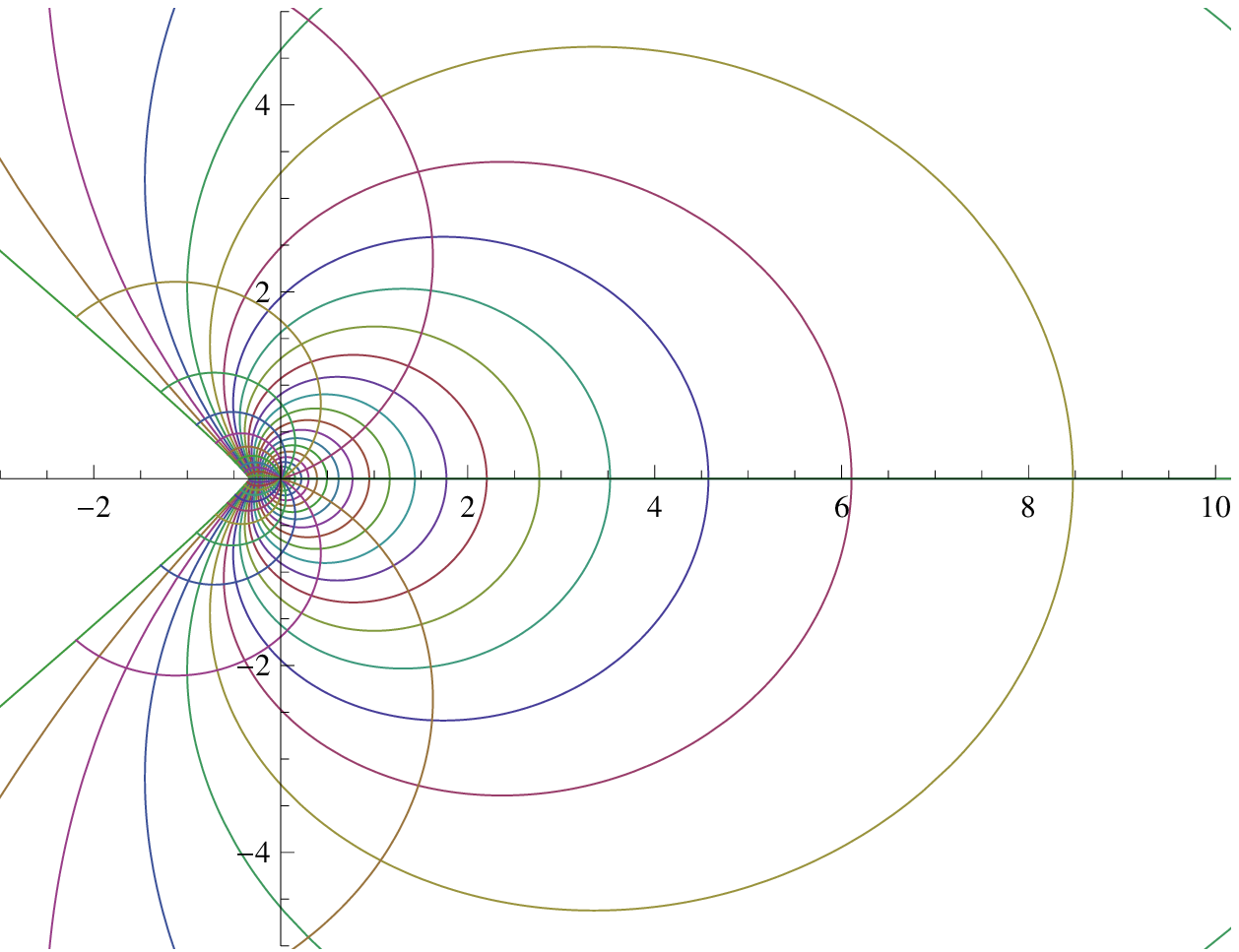}
\hspace{1cm}
\includegraphics[height=6.0cm, width=5.5cm, scale=1]{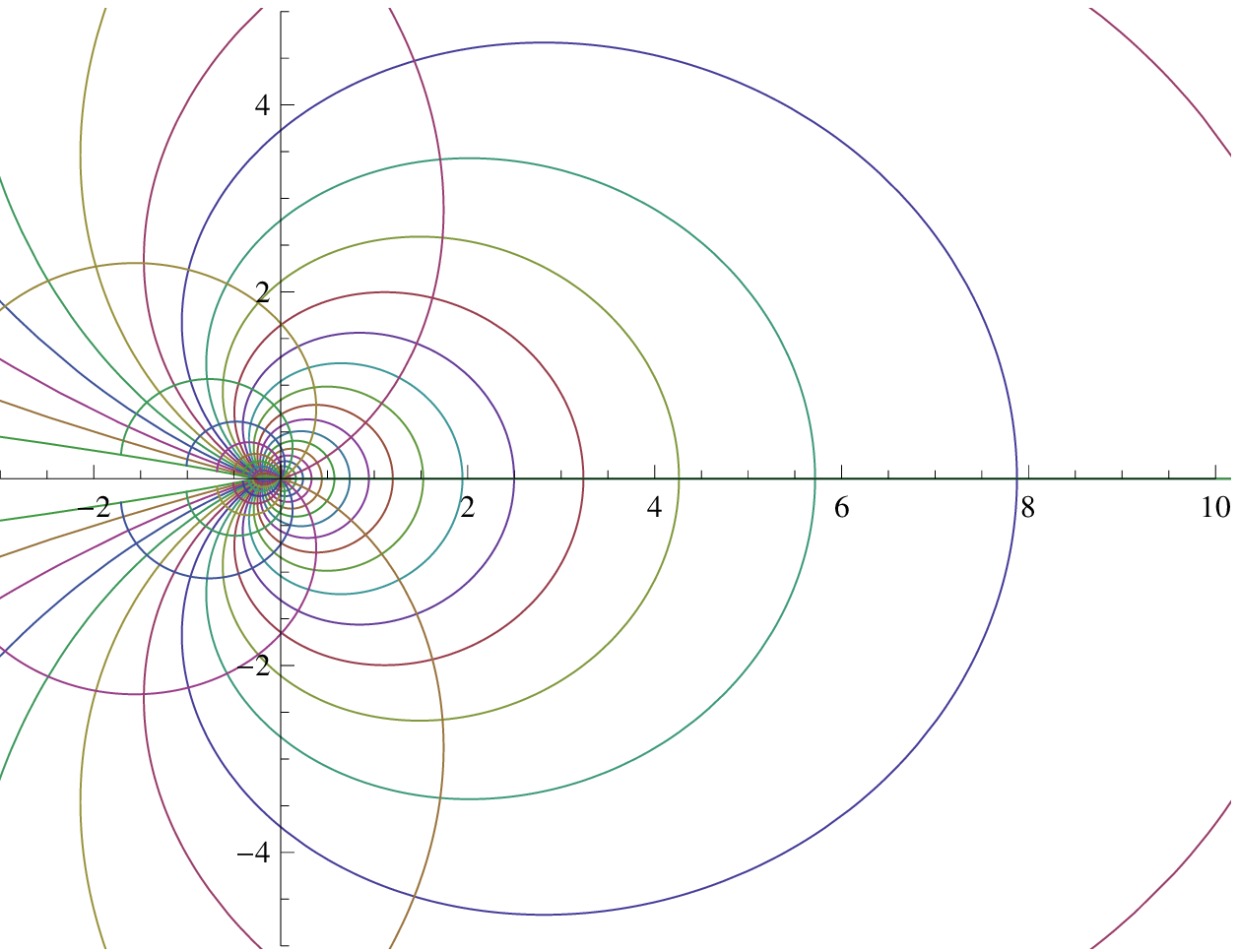}
\end{center}
(g) $\delta = 0.1 $, $K= 3$ \hspace{5cm} (h) $\delta = 0.1$, $K =3.85 $
\caption{The images of unit disk $\protect \ID$ under $f_{K,\delta}(z)=h(z)+\overline{g(z)}$ for certain values of $\delta$ and $K$
\label{fig2}}
\end{figure}

Images of $\mathbb{D}$ under the close-to-convex mappings $f_{K,\delta}(z)=h(z)+\overline{g(z)}$ for certain values of $\delta$ and $K$
with $2\leq K\leq 4-\delta $, where
$$ h(z)=\frac{1}{K}\left[ \left(\frac{1+z}{1-z}\right)^{K/2} - 1 \right] ~\mbox{ and }~ g'(z)=\sin\left(\frac{\delta \pi}{4}\right)zh'(z) ~\mbox{for $z\in\ID$},
$$
are drawn in Figure \ref{fig2}(a)--(h) using mathematica as plots of the images of equally spaced
radial segments and concentric circles of the unit disk.





Finally, we consider the class $ \mathcal{G}$ of functions $h\in {\mathcal A}$ such that
\be\label{cl-G}
{\rm Re } \left( 1+\frac{zh''(z)}{h'(z)}\right) < \frac{3}{2} ~\mbox{ for $z\in \ID$}.
\ee
Functions in $\mathcal{G}$ are known to be starlike in $\ID$. This class has been discussed recently, see for
example \cite{OPW13-SMJ} and the references therein. For this class we prove the following general result.

\bthm\label{AP-Theorem5}
Suppose that $h\in\mathcal{G}$ and satisfies the condition $ g'(z)=\omega (z) h'(z)$ in $\ID$, where
$\omega :\,\ID \rightarrow \ID$ is analytic, $\omega (0)=0$ and $W(z)=z(1+\omega (z))$ is starlike in $\ID$.
Then the harmonic mapping $f=h+\overline{g}$ is close-to-convex and univalent in $\ID$.
\ethm

\bcor\label{AP-Cor1}
Let $h\in\mathcal{G}$ and $g$ be analytic in $\ID$ such that $g'(z)=\lambda z^nh'(z)$ for some $n\in \IN$
and $0<|\lambda|\leq 1/{(n+1)}$.
Then the harmonic mapping $f=h+\overline{g}$ is close-to-convex and univalent in $\ID$.
\ecor
\bpf
Set $\omega (z)=\lambda z^n$ for $z\in\ID$. Then $W(z)=z(1+\omega (z))=z+\lambda z^{n+1}$ is starlike in $\ID$ if and only if
$|\lambda|\leq 1/{(n+1)}$. Indeed we have
$$|W'(z)-1|=(n+1)|\lambda| \, |z|^n<1 ~\mbox{ for $z\in \ID$}
$$
and hence, $W$ is univalent in $\ID$. Further, since $0<|\lambda|\leq 1/(n+1)$, we see that
$$\left |\frac{zW'(z)}{W(z)}-1 \right | =\left |\frac{n\lambda z^n}{1+\lambda z^n}\right |< \frac{n|\lambda|}{1-|\lambda|}\leq 1
~\mbox{ for $z\in \ID$}
$$
which implies that the function $W$ is starlike in $\ID$. The desired conclusion follows from Theorem~\ref{AP-Theorem5}.
\epf

\beg
According to Corollary \ref{AP-Cor1}, it follows that if $h\in\mathcal{G}$ and $g$ is analytic in $\ID$ such that $g'(z)=\lambda zh'(z)$
for some $\lambda$ with $ |\lambda| \leq 1/2$, then the harmonic mapping $f=h+\overline{g}$ is close-to-convex
(univalent) in $\ID$. For instance, let $h_1(z)=z-z^2/2$ and $g_1(z)=\lambda \left (\frac{z^2}{2}-\frac{z^3}{3}\right )$.
Then we see that $f_1=h_1+\overline{g_1}$ is clearly locally univalent in $\ID$ for each $|\lambda|< 1$. Also,
$$ 1+\frac{zh_1''(z)}{h_1'(z)}= \frac{1-2z}{1-z} ~\mbox{ for $z\in \ID$}
$$
and, since $w=(1-2z)/(1-z)$ maps $\ID$ onto the half-plane ${\rm Re}\, w <3/2$, by Corollary \ref{AP-Cor1}, it follows that
$f_1=h_1+\overline{g_1}$ is close-to-convex in $\ID$ for each $\lambda$ with $|\lambda | \leq 1/2$.
\eeg

We conjecture that Corollary \ref{AP-Cor1} is sharp in the sense that the bound on $\lambda$ cannot be improved.



\begin{figure}
\begin{center}
\includegraphics[height=5.5cm, width=5.5cm, scale=1]{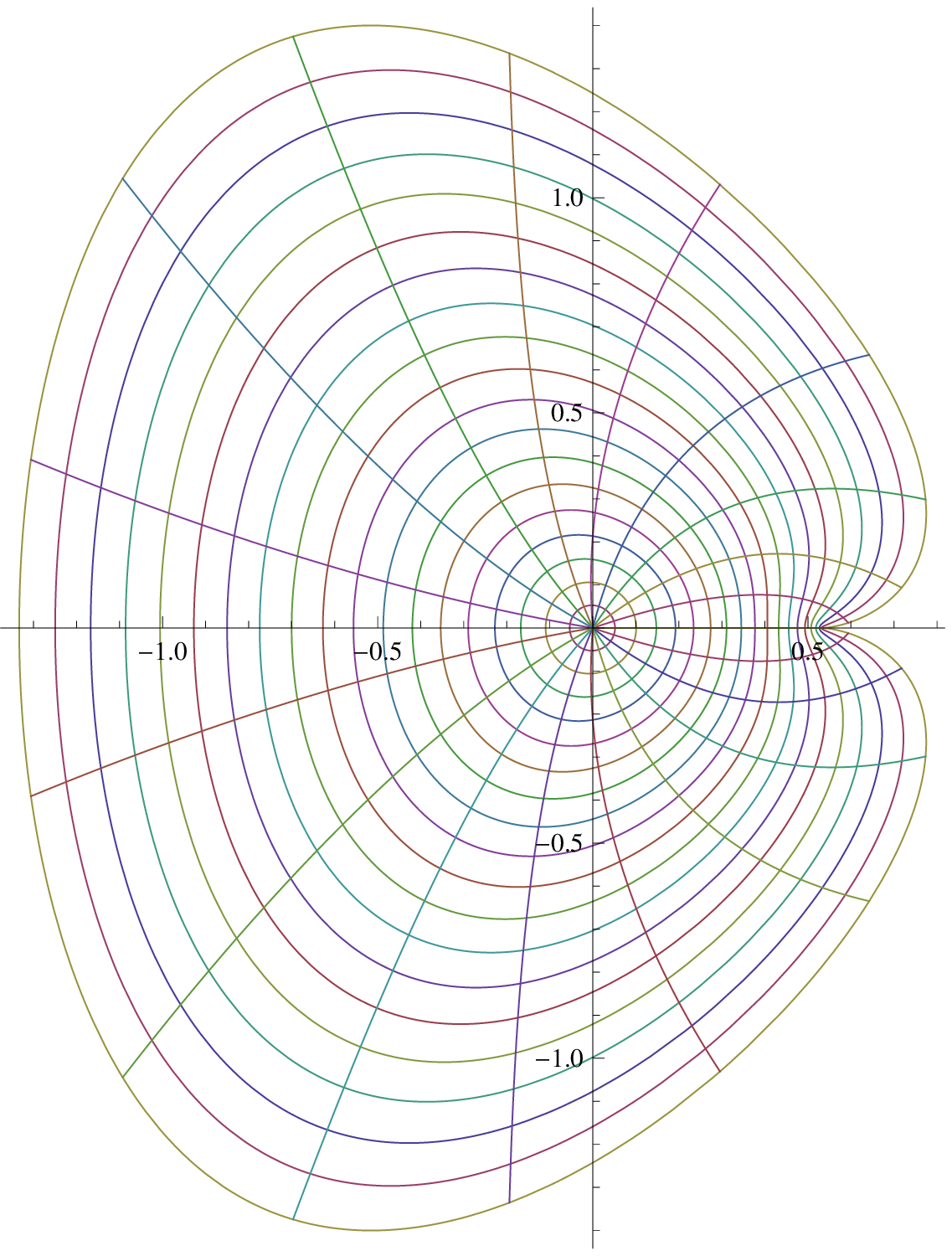}
\hspace{1cm}
\includegraphics[height=5.5cm, width=5.5cm, scale=1]{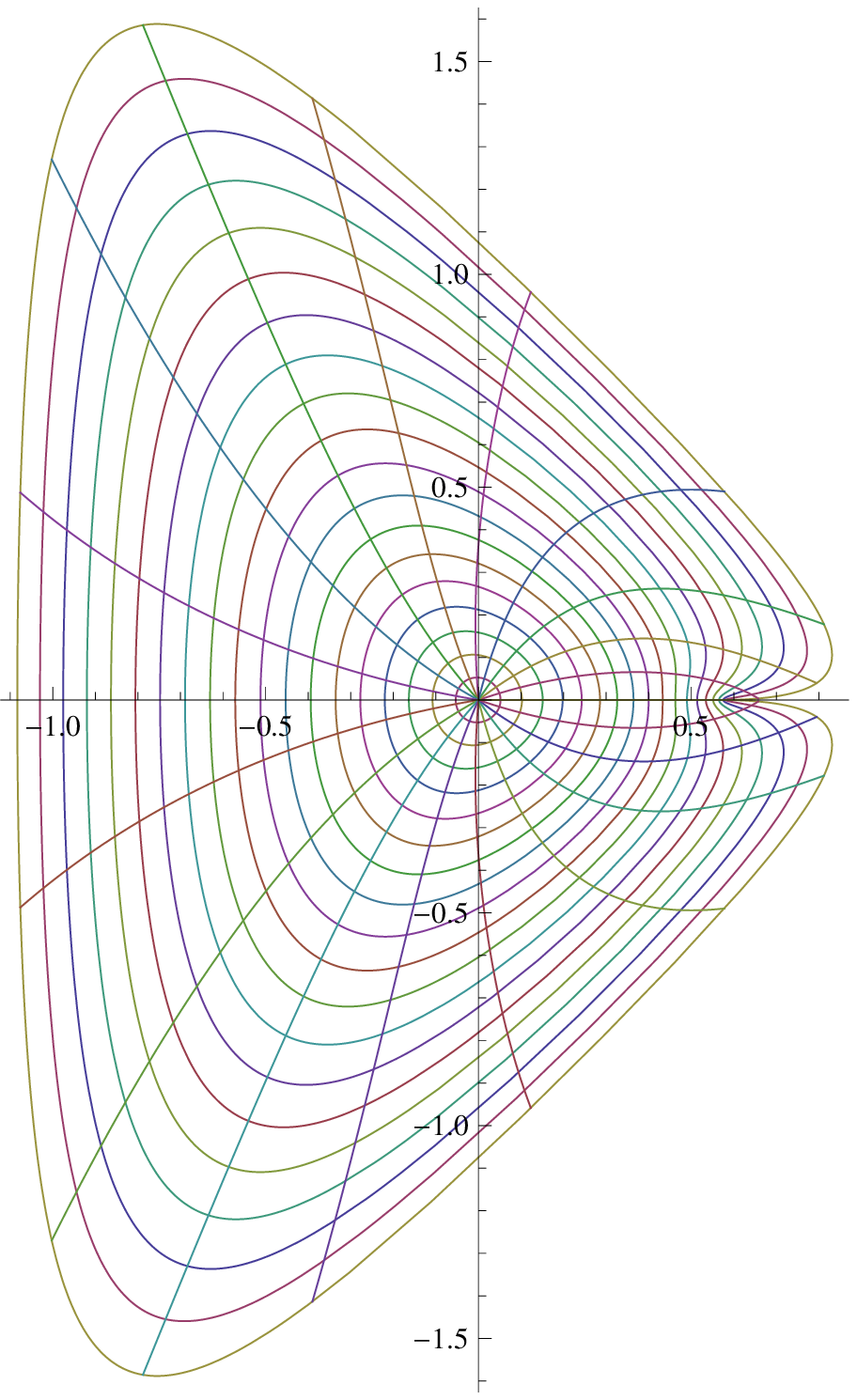}
\end{center}
{(a) $\lambda=\frac{1}{5}$ \hspace{5cm} (b) $\lambda=\frac{1}{2}$}


\begin{center}
\includegraphics[height=5.5cm, width=5.5cm, scale=1]{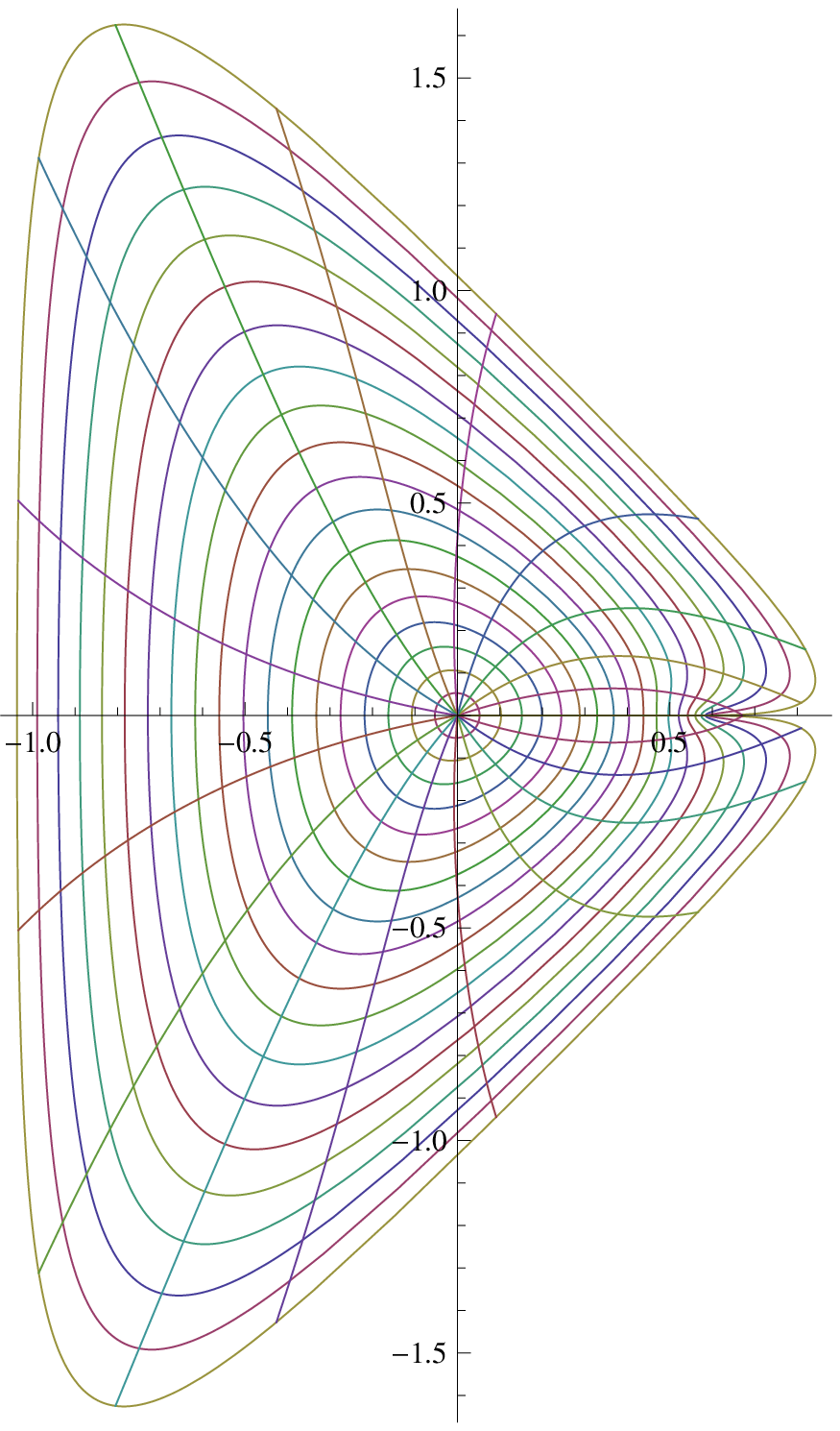}
\hspace{1cm}
\includegraphics[height=5.5cm, width=5.5cm, scale=1]{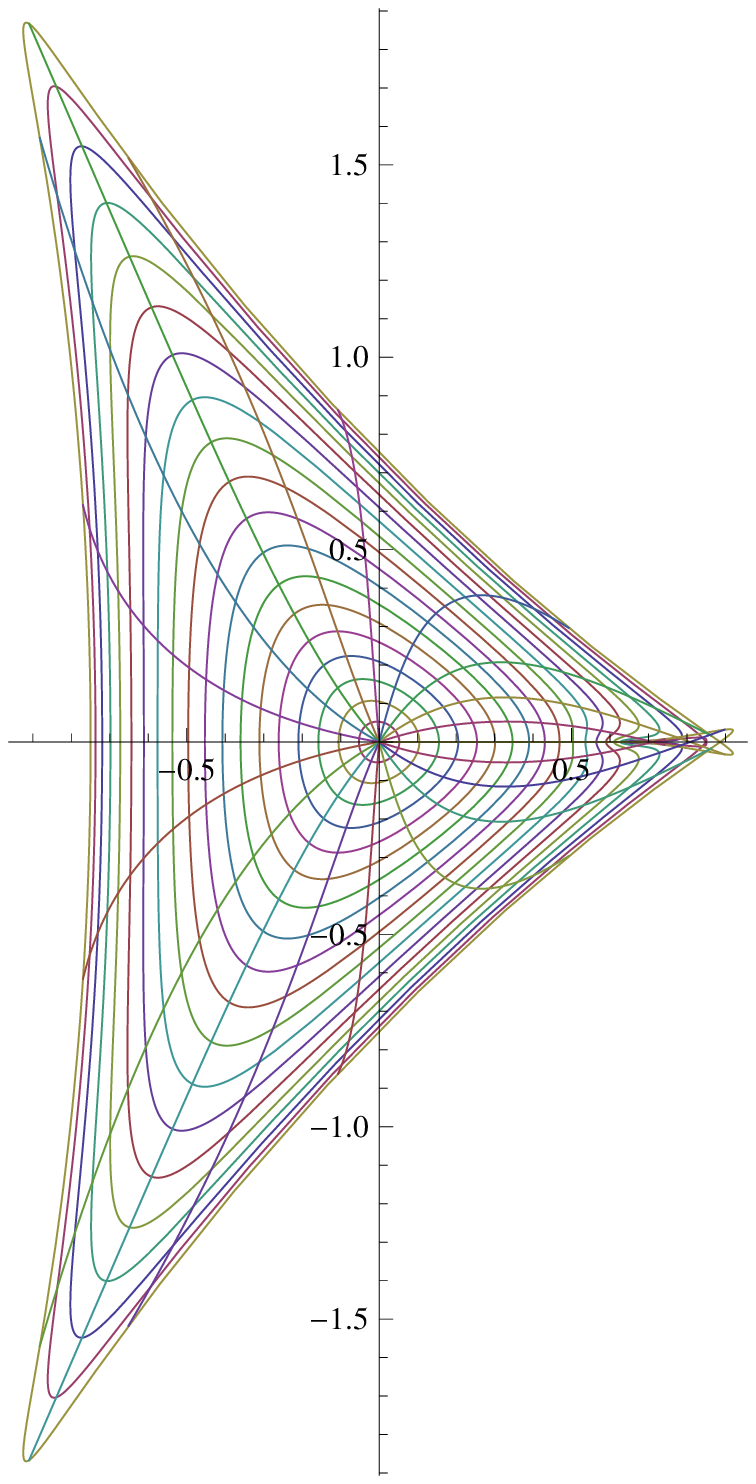}
\end{center}
(c) $\lambda=\frac{5}{9}$ \hspace{5cm} (d) $\lambda=\frac{9}{10}$
\caption{Image of $\mathbb{D}$ under $f(z)=z-\frac{z^2}{2}+\overline{\lambda \left (\frac{z^2}{2} -\frac{z^3}{3} \right )}$\label{fig1}}
\end{figure}
The image of unit disk under the function $f(z)=z-\frac{z^2}{2}+\overline{\lambda \left (\frac{z^2}{2} -\frac{z^3}{3} \right )}$ for
the values of $\lambda=1/5,1/2,5/9,9/10$ are drawn in Figure \ref{fig1}(a)--(d) using mathematica as plots of the images of equally spaced
radial segments and concentric circles of the unit disk. Closer examination of Figure \ref{fig1}(c)--(d) shows that the functions in these two
cases are not univalent in $\ID$.

\section{Proofs of Main Theorems}\label{sec-2}

\subsection{The class $\mathcal{F}$}
Let $ h\in \mathcal{K}(\beta)$ for some $-1/2\leq \beta <1$. 
Then
$$1+\frac{zh''(z)}{h'(z)}\prec p(z)=\frac{1+(1-2\beta )z}{1-z} ~\mbox{ for $z\in \ID$}.
$$
Here $\prec$ denotes the usual subordination (see \cite{Dur,Pomm}), and note that ${\rm Re}\,p(z)>\beta $ in $\ID$.
Thus, by the Hergtlotz representation for analytic functions with positive real part in the unit disk, it
follows that
$$\frac{zh''(z)}{h' (z)} =2(1-\beta) \int_{\partial \ID}\frac{\overline{x}z}{1-\overline{x}z }\,d\mu (x) ~\mbox{ for $z\in \ID$},
$$
where $\mu $ is a probability measure on $\partial \ID$ so that $\int_{\partial \ID} d\mu (x)=1$.
Therefore,
$$\log h'(z)=-2(1-\beta) \int_{\partial \ID} \log (1-\overline{x}z )\,d\mu (x) ~\mbox{ for $z\in \ID$}
$$
and thus, we have the sequence of functions $\{h_n(z)\}$ analytic in $\ID$,
\be\label{APW-eq3}
h_n'(z)= \prod _{k=1}^n(1-\overline{x}_{k}z)^{-2(1-\beta)t_k} ~\mbox{ where $|x_k|=1$, $0\leq t _{k}\leq 1$, $\sum\limits_{k=1}^{n}t _{k}=1$,}
\ee
which is dense in the family $ \mathcal{K}(\beta).$ The representation \eqref{APW-eq3} and Lemma \Ref{LemA} are two ingredients in the proof of
Theorem \Ref{BL-Theorem} and a similar approach helps to prove several other results.

\bpf[Proof of Theorem \Ref{BL-Theorem}]
Let $f=h+\overline{g}$, where $h\in {\mathcal K}(\beta)$ for some $\beta \in [-1/2, 0]$, and $ g'(z)=\omega (z) h'(z)$ for $z\in \ID$.
Then, it suffices to consider $h'(z)$ of the form
\be\label{eqh-1}
h'(z)= \prod _{k=1}^n(1-\overline{x}_{k}z)^{-2(1-\beta )t_k} ~\mbox{ for $z\in \ID$}.
\ee
Set
$$S(z)=\frac{z}{\prod _{k=1}^n(1-\overline{x}_{k}z)^{2t_k}} ~\mbox{ for $z\in \ID$}.
$$
Then $S$ is starlike in $\ID$.

\noindent\textbf{Case (1):} $\beta =-1/2$ and $\omega (z)=e^{i\theta}z$ for $z\in \ID$.

In this case, $ g'(z)=e^{i\theta}zh'(z)$ and thus, we may rewrite $h'(z)$ in the form
$$h'(z)= \prod _{k=1}^n(1-\overline{x}_{k}z)^{-3t_k}= \left (\prod _{k=1}^n\frac{1}{(1-\overline{x}_{k}z)^{t_k}} \right ) \frac{S(z)}{z} ~\mbox{ for $z\in \ID$},
$$
where $|x_k|=1$, $0\leq t _{k}\leq 1$ for $k=1,\ldots, n$, $\sum\limits_{k=1}^{n}t _{k}=1$.

Now, for each $|\epsilon |=1$ and $|\lambda |=1$, we consider the function
$$A(z)=\frac{z(h'(z)+\epsilon \overline{\lambda} g'(z))}{S(z)}
=\frac{1+\epsilon \overline{\lambda} e^{i\theta}z}{\prod _{k=1}^n(1-\overline{x}_{k}z)^{t_k}}
=\prod _{k=1}^n \left ( \phi_k(z) \right )^{t_k} ~\mbox{ for $z\in \ID$}, 
$$
where
$$ \phi_k(z)=\frac{1+\epsilon\overline{\lambda}e^{i\theta} z}{1-\overline{x}_{k}z} ~\mbox{ for $z\in \ID$}.
$$
Then the function $A(z)$ has the property that ${\rm Re } \left ( e^{i\gamma}A(z)\right )>0$ in $\ID$ for some $\gamma$.
In fact for each $k$, the function $\phi_k(z)$ maps the unit disk
$\ID$ onto a half plane so that ${\rm Re }\left (e^{i\theta _k}\phi_k(z)\right )>0$ for some $\theta _k$.
Setting $\gamma =\sum_{k=1}^nt_k\theta_k$, we have
$$|\arg \big (e^{i\gamma}A(z)\big )|=\left |\arg \prod _{k=1}^n \left ( e^{i\theta _k} \phi_k(z) \right ) ^{t_k}
\right | \leq \sum_{k=1}^nt_k \left |\arg \big (e^{i\theta _k}\phi_k(z)\big )\right | <\frac{\pi}{2}\sum_{k=1}^nt_k = \frac{\pi}{2}
$$
and hence, ${\rm Re } \left ( e^{i\gamma}A(z)\right )>0$ in $\ID$.
It follows that $F(z)=h(z)+\epsilon \overline{\lambda} g(z)$ is close-to-convex in $\ID$ for each $|\epsilon |=1$ and $|\lambda |=1$.
Hence, by Lemma \Ref{LemA}, functions $h+\lambda\overline{g}$ are close-to-convex in $\ID$,
for each $|\lambda |=1$. In particular, $f=h+\overline{g}$ is close-to-convex in $\ID$.

\noindent\textbf{Case (2):} Let $\beta =0$ and $|\omega (z)|<1$ for $z\in \ID$. In this case, using the relation 
$ g'(z)=\omega (z) h'(z)$ for $z\in\ID$, it follows easily that the function
$$\frac{z(h'(z)+\epsilon g'(z))}{S(z)}= \frac{zh'(z)(1+\epsilon \omega (z))}{S(z)} = 1+\epsilon \omega (z) 
$$
has positive real part in $\ID$ and hence, by Lemma \Ref{LemA}, the $h+\overline{g}$ is close-to-convex in $\ID$.

Next we assume that $\beta \in (-1/2, 0)$ and $|\omega (z)|< c=\cos (\beta \pi)$ for $z\in \ID$. In this case, using \eqref{eqh-1} we need to
consider the function
$$A(z)=\frac{z(h'(z)+\epsilon \overline{\lambda} g'(z))}{S(z)} =(1+\epsilon \overline{\lambda}\omega (z)) \prod _{k=1}^n(1-\overline{x}_{k}z)^{2\beta t_k}
$$
for each $|\epsilon |=1$ and $|\lambda |=1$. Again, it suffices to show that
${\rm Re } \left ( e^{i\gamma}A(z)\right )>0$ in $\ID$ for some $\gamma$.
In fact each fractional transformation of the form $w=\psi_k(z)= 1/(1-\overline{x}_{k}z)$ maps the unit disk
$\ID$ onto the half plane ${\rm Re}\, w>0$ so that
$$\left|\arg \big (1-\overline{x}_{k}z\big )^{2\beta t_k} \right |\leq \pi |\beta|t_k.
$$
and, since $|\arg (1+\epsilon \overline{\lambda}\omega (z))|<\arcsin c$ for $z\in\ID$, it follows that
$$|\arg A(z)|< \arcsin c + \pi |\beta|\sum_{k=1}^nt_k
= \arcsin c + \pi |\beta| =\pi/2.
$$
It follows that $F(z)=h(z)+\epsilon \overline{\lambda}g(z)$ is close-to-convex in $\ID$ for each $|\epsilon |=1$ and $|\lambda|=1$
and the desired conclusion follows from Lemma \Ref{LemA}.
\epf

\bpf[Proof of Theorem \ref{AP-Theorem}]
We can assume that $h$ belongs to the dense set of $ \mathcal{K}(-1/2)$, so in view of 
\eqref{APW-eq3} we have
$$h'(z)=\frac{1}{\prod _{k=1}^n(1-e^{-i\theta _k}z)^{3t _{k}}} ~\mbox{ for $z\in \ID$}, 
$$
where $\theta _k\in [0,2\pi]$, $0\leq t _{k}\leq 1$ and $\sum\limits_{k=1}^{n}t _{k}=1$. Using the last relation we find that
\beqq
\frac{1}{2\pi}\int_{0}^{2\pi}\frac{d\theta}{|h'(re^{i\theta})|^{2\beta }}
&= &\frac{1}{2\pi}\int_{0}^{2\pi} \prod _{k=1}^n\left |1-re^{i(\theta -\theta _k)}\right |^{6\beta t _{k}} \,d\theta\\
&\leq &\frac{1}{2\pi}\int_{0}^{2\pi} \sum _{k=1}^nt _{k}\left |1-re^{i(\theta -\theta _k)}\right |^{6\beta } \,d\theta\\
&=&\frac{1}{2\pi}\sum _{k=1}^nt _{k}\int_{0}^{2\pi} \left |1-re^{i(\theta -\theta _k)}\right |^{6\beta } \,d\theta\\
&=& \frac{1}{2\pi}\int_{0}^{2\pi} \left |1-re^{i\theta }\right |^{6\beta } \,d\theta\\
&\leq & 1+ \sum_{k=1}^\infty \left|{3\beta \choose k}\right |^2\\
&= &\frac{1}{2\pi}\int_{0}^{2\pi} \left |1-e^{i\theta }\right |^{6\beta } \,d\theta\\
&=& \frac{2^{3\beta}}{2\pi}\int_{0}^{2\pi} (1-\cos \theta )^{3\beta } \,d\theta\\
&=& \frac{2^{6\beta }}{\pi}\int_{0}^{\pi} \sin^{6\beta} (\theta) \, d\theta\\
&=& \frac{2^{6\beta +1}}{\pi}\int_{0}^{\pi/2} \sin^{6\beta} (\theta) \, d\theta\\
&=& \frac{2^{6\beta }}{\pi}{\bf B}\Big (\frac{6\beta +1}{2},\frac{1}{2}\Big ).
\eeqq
The desired conclusion follows. For the proof of sharpness part, we consider the function $h_0$ defined by
\be\label{APW-eq4}
h_0(z)=\frac{z-z^2/2}{(1-z)^2} 
=\frac{1}{2}\left (\frac{z}{1-z} +\frac{z}{(1-z)^2}\right ).
\ee
The function $h_0$ and its rotations belong to $\mathcal{F}$. A computation shows that $h_0'(z)=(1-z)^{-3}$
and the rest of the sharpness part follows easily.
\epf

\br As an alternate approach to the proof of Theorem \ref{AP-Theorem}, we may begin with $h\in \mathcal{F}$.
Then one has (see for instance \cite{samy-hiroshi-swadesh}) $h'(z)\prec (1-z)^{-3}$
and, since $h'(z)\neq 0$ in $\ID$, it follows that
$$\frac{1}{h'(z)} \prec (1-z)^3, ~z\in \ID.
$$
Thus (as in the proof of Theorem 1 in \cite{PoWi2014}) we see that
\beqq
\frac{1}{2\pi}\int_{0}^{2\pi}\frac{d\theta}{|h'(re^{i\theta})|^{2\beta }}
\leq \frac{1}{2\pi}\int_{0}^{2\pi} \left |1-re^{i\theta }\right |^{6\beta } \,d\theta
\eeqq
and the rest of the proof is as above. The desired result follows.
\er

%
%

\bpf[Proof of Theorem \ref{Thm-Alex2}]
By the hypothesis, there exists an analytic function $\omega:\,\ID \rightarrow \ID$ such that $g'(z)=\omega(z)h'(z)$. As $\omega (0)=g'(0)\in \ID$,
it follows that $|G'(0)|<|H'(0)|=1$. Moreover, since $h$ is starlike in $\ID$, $H$ is convex. Thus, according to Theorem \Ref{BL-Theorem}(2) with $\beta =0$,
it suffices to
show that $F_{\lambda}$ is sense-preserving in $\ID$. Indeed, since $ |g'(z)/h'(z)|<1$ in $\ID$, we obtain from a well-known result
of Robinson \cite[p.~30]{Robin1947} (see also \cite[Corollary 3.1]{MilMoc-1993}) that
$$|g(z)/h(z)|=|G'(z)/H'(z)| <1 ~\mbox{ in $\ID$}
$$
(and at the origin this is treated as the obvious limiting case). Thus, $F_{\lambda} =H+\lambda\overline{G}$ is
sense-preserving and harmonic in $\ID$. The desired conclusion follows from Theorem \Ref{BL-Theorem}(2).
\epf

\subsection{The class $CO(\alpha)$ of concave univalent functions}
We now consider normalized functions $h$ analytic in $\ID$ and
map $\ID$ conformally onto a domain whose complement with respect
to $\mathbb{C}$ is convex and that satisfy the normalization $h(1) = \infty$.
Furthermore, we impose on these functions the condition that
the opening angle of $h(\ID)$ at $\infty$ is less than or equal to
$\pi\alpha$, $\alpha\in (1,2]$. We will denote the family of such functions by $CO(\alpha)$
and call it as the class of concave univalent functions.
We note that for $h\in CO(\alpha)$, $\alpha\in (1,2]$, the closed set $\IC\backslash h(\ID)$
is convex and unbounded. Also, we observe that $CO(2)$ contains the classes $CO(\alpha)$,
$\alpha\in (1,2]$.

For the proof of Theorem \ref{AP-Theorem3}, we need the following result due to Avkhadiev and Wirths \cite{Avk-Wir-05}.

\begin{Thm}(\cite[Theorem 1]{Avk-Wir-05}) \label{TheoB}
The set of functions $h\in CO(\alpha ),$ with
\be\label{AP-eq1}
h'(z)=\frac{\prod\limits_{k=1}^{n}\left( 1-e^{it_k}z\right) ^{\beta
_{k}}}{( 1-z)^{\alpha +1}},
\ee
where $0 <t_1<t_2<\cdots < t _{n}<2\pi $,
$0<\beta _{k}\leq 1$ for $k=1,\ldots, n$, with
$\sum\limits_{k=1}^{n}\beta _{k}=\alpha -1,$ is dense in $CO(\alpha ).$
\end{Thm}

\bpf[Proof of Theorem \ref{AP-Theorem3}]
Let $h\in CO(\alpha )$. Then, according to Theorem \Ref{TheoB}, it suffices to prove the theorem for $h$ of the form \eqref{AP-eq1},
where $0 <t_1<t_2<\cdots < t _{n}<2\pi $, $0<\beta _{k}\leq 1$ for $k=1,\ldots, n$ with
$\sum\limits_{k=1}^{n}\beta _{k}=\alpha -1.$

Now, we set $F=h+\epsilon g$, where $|\epsilon| =1$.
Then, because $ g'(z)= \omega (z)h'(z)$ in $\ID$, the above representation on $h'(z)$ gives that
\[F'(z)= 
(1+\epsilon\omega (z))h'(z)=\frac{(1 +\epsilon\omega (z))}{(1-z)^{2}}
\prod\limits_{k=1}^{n}\big ( \phi_k(z)\big )^{\beta _{k}}
\]
where
$$ \phi_k(z)=\frac{1-e^{it_k}z}{1-z}.
$$
With $k(z)=\frac{\overline{c}z}{(1-z)^2}$ with $|c|=1$, it follows that
\[\frac{zF'(z)}{k(z)} 
=c(1 +\epsilon \omega (z))\prod\limits_{k=1}^{n} \big ( \phi_k(z)\big )^{\beta _{k}}.
\]
We observe that each $ \big ( \phi_k(z)\big )^{\beta _{k}}$ forms a
wedge at the origin with angle of measure $\beta _{k}\pi/2 $ and containing the point $1$.
Hence the product with $c$ make angles of total less than $(\alpha -1)\pi /2$.
Next, we note by hypothesis that $|\omega (z)|<c=\sin (\frac{2-\alpha }{2})\pi $ for $z\in\ID$,
and thus, we deduce that $|\arg (1+\epsilon \omega(z) )|<\arcsin c \, =(2-\alpha )\frac{\pi }{2}$.
Thus,
$$\left |\arg \left (c(1+\epsilon \omega(z))\prod\limits_{k=1}^{n} \big ( \phi_k(z)\big )^{\beta _{k}} \right )\right |<(2-\alpha )\frac{\pi }{2} +
(\alpha -1)\frac{\pi }{2} =\frac{\pi }{2}.
$$
Observe that the existence of an unimodular complex constant $c$ is guaranteed as in the proof of Theorem \Ref{BL-Theorem}. Therefore,
$$ {\rm Re}\, \left (c (1-z)^2F'(z)\right ) = {\rm Re}\, \left (\frac{zF'(z)}{k(z)}\right ) >0 \mbox{ for $z\in\ID$}
$$
and hence, for each $\epsilon$ with $|\epsilon| =1$, the analytic function $F=h+\epsilon g$
is close-to-convex in $\ID$. The desired conclusion follows from Lemma \Ref{LemA}.
\epf

\subsection{The class ${\mathcal V}_K$ of functions of bounded boundary rotation}

For the proof of Theorem \ref{AP-Theorem4}, we need some preparation. We begin to recall
 the familiar representation obtained by Paatero \cite{Paa1} for functions $h\in {\mathcal V}_K$:
\begin{equation}\label{sec2-eq1a}
 h'(z)=\exp\left(-\int_0^{2\pi}\log(1-ze^{-it})\,d\mu(t)\right),
\end{equation}
where $\mu(t)$ is a real valued function of bounded variation on $[0,2\pi]$ with
\be\label{sec2-eq1aa}
\int_0^{2\pi} d\mu (t)=2 ~\mbox{ and }~ \int_0^{2\pi}|d\mu (t)|\le K .
\ee
It is well-known that ${\mathcal V}_2$ coincides with the class of normalized convex univalent functions and that
for $2\leq K \leq 4$, all members of ${\mathcal V}_K$ are univalent in $\ID$ (see \cite{Paa1}).
However, Pinchuk \cite{Pin1} strengthen this result by showing that for $2\leq K \leq 4$, the classes ${\mathcal V}_K$
consist of all close-to-convex functions. This fact also follows from our Theorem \ref{AP-Theorem4}.
However, each of the classes ${\mathcal V}_K$ with $K>4$ contains non-univalent functions.
An extremal function belonging to this class is
\begin{equation}\label{sec2-eq2}
g_K(z)=\frac{1}{K}\left[\left(\frac{1+z}{1-z}\right)^{K/2}-1\right] .
\end{equation}
It had been shown by Pinchuk \cite[Theorem~6.2]{Pin2} that the image of the unit disk $\ID$ under a ${\mathcal V}_K$ function
contains the disk of radius $1/K$ centered at the origin, and the functions of the class ${\mathcal V}_K$ are
continuous in $\overline\ID$ with the exception of at most $[K/2+1]$ points on the unit circle $\partial \ID$.
Moreover, it is known that $h\in {\mathcal V}_K$ if and only if
$$1+\frac{zh''(z)}{h'(z)}=\left (\frac{K}{4}+\frac{1}{2} \right )p_1(z) -\left (\frac{K}{4}-\frac{1}{2} \right )p_2(z)
$$
for some $p_1,p_2\in {\mathcal P}$, where ${\mathcal P}$ denotes the class of functions
$p$ analytic in $\ID$ such that $p(0)=1$ and ${\rm Re}\, p(z)>0$ in $\ID$. Thus, it follows that
\begin{equation*}
h'(z)=\exp \left ( -\int_{|x|=1} \log (1-xz)\left( d\mu _{1}(x)-d\mu _{2}(x)\right)\right )
\end{equation*}
where
$$\int _{|x|=1}\left( d\mu _{1}(x)-d\mu _{2}(x)\right) =2, ~~\int _{|x|=1} d\mu _{1}(x)=\frac{K}{2}+1, ~\mbox{ and }~
\int_{|x|=1} d\mu _{2}(x)=\frac{K}{2}-1.
$$
Consequently, we easily have the following

\begin{lem}\label{lem1a}
If $h\in{\mathcal V}_K$, then there exists a sequence of functions $\{h_n(z)\}$ analytic in $\ID$
such that
\be\label{APW-eq2}
h_{n}'(z)=\frac{\prod\limits_{k=1}^{n}(1-\overline{x}_{k}z)^{ \alpha _{k}}}{\prod\limits_{k=1}^{n}(1-\overline{y}_{k}z)^{ \beta _{k}} },
\ee
where $|x_{k}|=1$, $|y_{k}|=1$, $0\leq \alpha _{k},\,\beta _{k} \leq 1$ with
\be\label{APW-eq2e}
\sum\limits_{k=1}^{n}\alpha _{k}=\frac{K}{2}-1 ~\mbox{ and }~\sum\limits_{k=1}^{n}\beta _{k}=\frac{K}{2}+1,
\ee
and $\{h_n\}$ converges uniformly on compact subsets of $\ID$. That is, $\{h_n(z)\}$ is dense in the family
${\mathcal V}_K.$
\end{lem}

\bpf[Proof of Theorem \ref{AP-Theorem4}]
Set $F=h+\epsilon g$, where $|\epsilon|=1$ and $h\in{\mathcal V}_K$. In view of Lemma \ref{lem1a}, it suffices to choose $h\in {\mathcal V}_{K}$ so that
$$h'(z)=\frac{\prod\limits_{k=1}^{n}(1-\overline{x}_{k}z)^{ \alpha _{k}}}{\prod\limits_{k=1}^{n}(1-\overline{y}_{k}z)^{ \beta _{k}} },
$$
where $|x_{k}|=1$, $|y_{k}|=1$, $0\leq \alpha _{k},\,\beta _{k} \leq 1$ satisfying the conditions \eqref{APW-eq2e}. It is convenient to rewrite the
last expression as
$$h'(z)= \prod\limits_{k=1}^{n} \left (\frac{1-\overline{x}_kz}{1-\overline{y}_kz}\right )^{ \alpha _k} \cdot
\prod\limits_{k=1}^{n}(1-\overline{y}_{k}z)^{ t_k}, \quad t_k=\alpha _{k}- \beta _{k}.
$$
Observe now that $\sum\limits_{k=1}^{n}t_k=2$ and thus, the function $S$ defined by
$$S(z)=\frac{\overline{c}z}{\prod \limits_{k=1}^{n}(1-\overline{y}_{k}z)^{t _k}} \quad (|c|=1)
$$
is starlike in $\ID$. Further,
$$ F'(z)=h'(z)+\epsilon g'(z)=(1+\epsilon \omega (z))h'(z)
$$
so that
$$\frac{zF'(z)}{S(z)}=c(1+\epsilon \omega (z))\prod\limits_{k=1}^{n} \left (\frac{1-\overline{x}_kz}{1-\overline{y}_kz}\right )^{ \alpha _k}
$$
where (by the hypothesis)
$$\sum\limits_{k=1}^{n} \alpha _{k}=\frac{K}{2}-1\leq 1-\frac{\delta }{2}.
$$
Note that $ |\arg (1+\epsilon\omega (z))|<\pi\delta /4$.
This observation shows that (with a suitably defined $c$ on $\partial \ID$)
$$ \left |\arg \left (\frac{zF'(z)}{S(z)}\right )\right |
< \frac{\pi}{2} \left (\sum\limits_{k=1}^{n} \alpha _{k}+ \frac{\delta }{2}\right )=
\frac{\pi}{2} \left ( \frac{K}{2}-1+ \frac{\delta }{2}\right ) \leq \frac{\pi}{2}
$$
and thus, the function $zF'(z)/S(z)$ has positive real part in $\ID$. It follows that
$F(z)=h(z)+\epsilon g(z)$ is close-to-convex in $\ID$
for each $|\epsilon |=1$ and hence, by Lemma \Ref{LemA}, the
harmonic function $ f=h+\overline{g}$ is close-to-convex in $\ID$.
\epf

\subsection{The class $\mathcal{G}$}

Now, we let $h\in\mathcal{G}$. Then \eqref{cl-G} holds. Clearly, \eqref{cl-G} can be written as
$$ 1+ \frac{zh''(z)}{h'(z)}\prec p(z)=\frac{1-2z}{1-z} ~\mbox{ for $z\in \ID$}
$$
and thus, by the Hergtlotz representation for analytic functions with positive real part in the unit disk, it
follows easily that
$$\frac{h''(z)}{h' (z)} =- \int_{\partial \ID}\frac{\overline{x}}{1-\overline{x}z }\,d\mu (x)~\mbox{ for $z\in \ID$},
$$
where $\mu $ is a probability measure on $\partial \ID$ so that $\int_{\partial \ID} d\mu (x)=1$. This
means that
$$ h'(z)= \exp \left (\int_{\partial \ID} \log (1-\overline{x}z )\,d\mu (x) \right )~\mbox{ for $z\in \ID$}.
$$
Thus, we have a sequence of functions $\{h_n(z)\}$ analytic in $\ID$ such that
\be\label{APW-eq3a}
h_n'(z)= \prod _{k=1}^n(1-\overline{x}_{k}z)^{\alpha _k}
\ee
where $|x_k|=1$, $0\leq \alpha _{k}\leq 1$ for $k=1,2,\ldots , n$, $\sum\limits_{k=1}^{n}\alpha _{k}=1$,
and $h_n\ra h$ uniformly on compact subsets of $\ID$. That is, $\{h_n(z)\}$ is dense in the family
${\mathcal G}.$ We observe that functions in ${\mathcal G}$ are bounded in $\ID$.

\bpf[Proof of Theorem \ref{AP-Theorem5}]
As in the proofs of previous theorems, we begin to set $F=h+\epsilon g$, where $|\epsilon|=1$ and $h\in {\mathcal G}.$
In view of the above discussion and \eqref{APW-eq3a}, it suffices to prove the theorem for
functions $h$ of the form
$$h'(z)= \prod _{k=1}^n(1-\overline{x}_{k}z)^{\alpha _k}
$$
where $|x_k|=1$, $0\leq \alpha _{k}\leq 1$ for $k=1,2,\ldots , n$ and $\sum\limits_{k=1}^{n}\alpha _{k}=1$.
Consequently, there exists a complex number $c$ with $|c|=1$ and such that
$$\frac{czF'(z)}{W(z)}=c\prod \limits_{k=1}^n(1-\overline{x}_{k}z)^{\alpha _k}
$$
has positive real part for $z\in\ID$, where $W$ defined by $W(z)=z+\epsilon z\omega (z)$ is starlike for each $|\epsilon|=1$
(by hypothesis). Thus, the harmonic function $ f=h+\overline{g}$ is close-to-convex in $\ID$ (by Lemma \Ref{LemA}).
\epf

\end{document}